\documentclass[12pt]{amsart}
\usepackage{amsfonts,amsmath,amssymb,comment,amscd,latexsym,}
\numberwithin{equation}{section}

\newtheorem{Lem}{Lemma}[section]
\newtheorem{Prop}[Lem]{Proposition}
\newtheorem{Cor}[Lem]{Corollary}
\newtheorem{Thm}[Lem]{Theorem}

\theoremstyle{definition}
\newtheorem{Def}[Lem]{Definition}

\renewcommand{\k}{{\Bbbk}}
\renewcommand{\o}{\otimes}
\newcommand{\e}{e_{\lambda}}
\newcommand\gen[1]{{\langle #1\rangle}}
\newcommand{\G}{{\widehat G}}
\newcommand{\D}{\mathcal{D}}
\newcommand{\pf}{\medskip\noindent{\sc Proof:\,}}
\newcommand{\1}{\otimes 1_{\lambda}}
\newcommand{\R}{R(\lambda)}
\newcommand{\Z}{Z(\lambda)}

\begin{document}

\title[Trace-Like Invariant]{Trace-Like Invariant for Representations of Nilpotent Liftings of Quantum Planes}
\author{Andrea Jedwab}
\address{University of Southern California\\Los Angeles, CA 90089}\email{jedwab@usc.edu} 
\author{Leonid Krop}
\address{DePaul University\\Chicago, IL 60614}\email{lkrop@condor.depaul.edu}

\begin{abstract}We derive a formula for the trace of the antipode on endomorphism algebras of simple self-dual modules of nilpotent liftings of quantum planes. We show that the trace is equal to the quantum dimension of the module up to a nonzero scalar depending on the simple module.
\end{abstract}
\date{7/30/09}
\maketitle

\section*{Introduction}\label{intro}

The classical Frobenius-Schur function (or indicator) $\nu_2(\chi)$ on the set $\{\chi\}$ of complex characters of a finite group $G$ \cite{F} has been extended to simple modules of a semisimple finite-dimensional Hopf algebra $H$ by Linchenko and Montgomery \cite{LM}. Let $S$ be the antipode of $H$, $V$ a simple self-dual $H$-module. Denote by $S_V$ the endomorphism induced by $S$ on $\text{End}(V)$. It was noted in \cite{LM} that in the semisimple case, as a consequence of the fact that $S^2=\text{id},\ \nu_2(V)=\displaystyle \frac{\text{Tr}(S_V)}{\text{dim}V}$ and $\text{Tr}(S_V)=\pm\text{dim}\ V,$ depending on whether the indicator is $1$ or $-1.$ In particular, $\text{Tr}(S_V)$ is an integer closely related to the Frobenius-Schur indicator. In general, we can still consider the map $S_V$ on $\text{End}(V)$ for self-dual modules of any finite dimensional Hopf algebra. Defined as $\mu(V)$ in \cite{J}, $\text{Tr}(S_V)$  was deemed a good subsitute of the indicator in the non-semisimple case. Indeed, at least for some Hopf algebras, a similar relation exists between $\text{Tr}(S_V)$ and the quantum dimension of $V$ with the coefficient $\pm 1$ modified by a root of unity.

In this note we work with Hopf algebras which are no longer semisimple. They belong to the class of algebras called nilpotent type liftings of a quantum plane \cite{CK}, a subfamily of Hopf algebras introduced and classified in \cite{AS1}. The defining relations of these algebras are encoded in a tuple called datum which involves as a parameter a certain root of unity $q$. Our main result is that for a self-dual module $V$, $\text{Tr}\ (S_V)=\rho(V)\displaystyle{(-1)^{\text{dim}V-1}}(\text{dim}\ V)_q$ where $\rho(V)$ is a root of $1$ depending on $V$ and, for a natural number $m$, $(m)_q$ denotes the $q$- Gaussian integer.

\section{Background Review}\label{prelims}

\subsection{Liftings of Quantum spaces}\label{liftings}

We fix some notation. Below $\k$ is a field of characteristic $0$ containing all roots of $1$ and $\k^{\bullet}=\k\setminus\{0\}$. We let $G$ stand for a finite abelian group, $\G:={\rm Hom}(G,\k^{\bullet})$ denote the dual group, and write $\k G$ for the group algebra of $G$ over $\k$. The order of $g\in G$ is denoted by $|g|$ and $|q|$ denotes the order of an element $q\in\k^{\bullet}$. For an $H$- module $V$ (algebra $A$) and a subset $X$ of $V$ (A) we denote by $\gen X$ the submodule (subalgebra) generated by $X$. The unsubscribed `$\otimes$' means `$\otimes_{\Bbbk}$'. For all $n,m\in\mathbb Z$ with $m\ge 0$, ${\binom nm}_q$ denotes the Gaussian $q$-binomial coefficient \cite{KC}, $(n)_q={\binom n1}_q$ and $(n)_q!=(1)_q\cdots (n)_q$.

Liftings of quantum linear spaces were constructed and completely described in \cite{AS1}. Their representation theory was studied in e.g. \cite{{AB1},{AB2},{CK},{KR}}. When the dimension of a quantum linear space is two, we call them liftings of quantum planes. In the terminology of \cite{{CK},{KR}} a nilpotent type linked lifting of a quantum plane is a Hopf algebra whose defining relations are determined by a tuple 
$$\D=(G, a,b,\chi,1,0)$$
called \emph{datum}, composed of a finite abelian group $G$, elements $a,b\in G$, and a character $\chi\in G$. The datum satisfies the conditions:
\begin{align}n&=|\chi(a)|>1\\
             \chi(a)&=\chi(b)\\
             ab&\ne 1\end{align}
Let $q=\chi(a^{-1})$ so that $n=|q|$. The Hopf algebra $H=H(\D)$ associated to $\D$ is generated by $G,$ $x$ and $y$ subject to the relations of $G$ and the following additional relations
\begin{align}x^n&=0=y^n\label{xypower}\\
              gx&=\chi(g)xg\;\text{for all}\;g\in G\label{char1}\\
              gy&=\chi^{-1}(g)yg\; \text{for all}\;g\in G\label{char2}\\
              xy &-qyx=ab-1\label{qcom}\end{align}
 The coalgebra structure of $H$ is given by
 \begin{align}\Delta(x)&=a\otimes x+x\otimes 1\\
                \Delta(y)&=b\otimes y+y\otimes 1\\
                \Delta(g)&=g\otimes g\;\text{for all}\;g\in G\\
                \epsilon(x)&=0=\epsilon(y)\\
                \epsilon(g)&=1\;\text{for all}\;g\in G,
                \end{align}             
and the mapping $S:H\to H$ defined on generators by $S(g)=g^{-1},g\in G, S(x)=-a^{-1}x, S(y)=-b^{-1}y$ extends to the antipode of $H$.

\subsection{Simple Modules}\label{simples}

We recall a standard theory of simple $H$- modules following \cite{KR}. Consider the sub-Hopf algebra  $Y=\gen{G,y}$ and set $Y'=\sum\k Gy^i$. Clearly $Y=\k G\oplus Y'$ and $Y'^n=0$, hence $Y'$ is the radical of $Y$. Thus every simple $Y$- module is a pull-back of a simple $\k G$- module. The latter is defined by a character $\lambda:G\to\k^{\bullet}$. Let $\k_{\lambda}$ denote the simple $Y$- module associated to $\lambda$ and set $1_{\lambda}=1_{\k}$. Then $\k_{\lambda}$ is described by the relations
\begin{align*}g.1_{\lambda}&=\lambda(g)1_{\lambda}\;\text{for all}\;g\in G\\
              y.1_{\lambda}&=0\end{align*}
Simple $H$-modules are obtained by way of certain induced modules, called Verma-type modules. To every $\lambda\in\G$ we associate a cyclic $H$- module $Z(\lambda)$ 
\begin{equation} Z(\lambda)=H\otimes_Y\k_{\lambda}\label{Verma}\end{equation}
Usefulness of $\Z$ comes from its universal property. Say $M$ is an $H$- module.  We say that $m$ is a $\lambda$- weight element if $g.m=\lambda(g)m$ for all $g\in G$. An element $m$ of $M$ is called Y- primitive (or just primitive when $Y$ is clear from the context) if $ym=0$, and $m$ is $\lambda$- primitive if it also has weight $\lambda$. Suppose $M$ is a cyclic $H$- module generated by a primitive  element $m_0$ of weight $\lambda$. Then there is a unique epimorphism $\Z\twoheadrightarrow M$ sending $1\otimes 1_{\lambda}$ to $m_0$. Now let $L$ be a simple $H$- module. Restricting $L$ to $Y$ we pick a simple nonzero $Y$- module in the $Y$- socle of $M$. Since simple $Y$- modules are $1$- dimensional, this submodule is $\k m$ for a primitive $m$ of certain weight $\lambda$. It follows that $L$ is an epimorphic image of $Z(\lambda)$. In fact, the weight $\lambda$ is uniquely determined by $L$, hence we write $L=L(\lambda)$. We summarize the above discussion in the next theorem and, for the reader's convenience, give a self-contained proof (see also \cite[3.5]{KR}).
\begin{Thm}\label{parametrization} The set $\{L(\lambda)|\lambda\in\G\}$ is the full set of representatives of isomorphism classes of simple $H$- modules.
\end{Thm}
\pf By \cite[5.2]{AS1} $H$ has a basis $\{x^igy^j|0\le i,j\le n-1, g\in G\}$. It follows that, by its construction, $Z(\lambda)$  has a basis $\{x^i\otimes 1_{\lambda}|0\le i\le n-1\}$. This is the \emph{standard} basis of $Z(\lambda)$. The equality $g.(x^i\otimes 1_{\lambda})=\lambda\chi^i(g)(x^i\otimes 1_{\lambda})$ means that $x^i\otimes 1_{\lambda}$ has weight $\lambda\chi^i$. Since $\chi(a)=q,\;|\chi|\ge n$, and therefore all $\lambda\chi^i$ are distinct. We conclude that $Z(\lambda)$ has a unique line of vectors of weight $\lambda$. It follows that a sum of proper submodules is proper. For, otherwise one of them would contain $1\1$ because a weight subspace of a given weight of a sum of modules is the sum of the weight subspaces of this weight of the summands. Let $\R$ denote the unique maximal submodule of $\Z$. Note that every weight vector is a monomial $x^j\1,j\ge 1$, hence a primitive weight vector generates a proper submodule, namely the span of all $x^k\1,k\ge j$. Therefore $L(\lambda):=\Z/\R$ doesn't have a primitive weight vector of a weight distinct from $\lambda$. Since an isomorphism carries $\lambda$- primitives into $\lambda$- primitives, $L(\lambda)\cong L(\mu)$ iff $\lambda=\mu$.\qed

We can say a great deal more about $L(\lambda)$ by computing $\R$. First we make
\begin{Def}\label{exponent} For a weight $\lambda\in\G$ we define a nonnegative integer $e=e(\lambda)$ by
$$\lambda(ab)=q^e\;\text{if}\;\lambda(ab)\in\gen q\;\text{and}\;e(\lambda)=n-1,\;\text{otherwise}.$$
\end{Def}
There $\gen q$ denotes the subgroup of $\k^{\bullet}$ generated by $q$.
\begin{Prop}\label{radical} $\R$ is a proper submodule iff $e(\lambda)\le n-2$. If so $\dim\ L(\lambda)=e+1$.
\end{Prop}
\pf Since every $H$- module contains a primitive weight element, so does $\R$. To find that element we need an identity for skew(braided)- commutators. This is \cite[1.2(2)]{KR}, which in view of the relation $yx-q^{-1}xy=(-q^{-1})(ab-1)$ requires changing $q$ by $q^{-1}$ and setting $\lambda=-q^{-1}$. Thus we have
\begin{equation}\label{simplecommutator}yx^s=q^{-s}x^sy-q^{-s}(s)_qx^{s-1}(q^{-(s-1)}ab-1).\end{equation}
Suppose $x^s\1,s\ge 1$ is primitive. That is to say
\begin{align*}0=y(x^s\1)=yx^s\1&=-q^{-s}(s)_qx^{s-1}(q^{-(s-1)}ab-1)\1\\
                               &=-q^{-s}(s)_q(q^{-(s-1)}\lambda(ab)-1)x^{s-1}\1\end{align*}
This equality implies $\lambda(ab)=q^{s-1}$ which says in turn that $\Z$ is reducible (i.e. $\R\ne 0$) iff $\lambda(ab)\in\{q^i|0\le i\le n-2\}$, as $1\le s\le n-1$. Now assume $\lambda(ab)=q^e$ with $e< n-1$. Then $x^{e+1}\1$ is the only primitive element within a scalar multiple distinct from $1\1$. Thus $\R=H(x^{e+1}\1)=\gen{x^j\1|j\ge e+1}$, whence $\dim\,L(\lambda)=e+1$ and the proof is complete.\qed 

Below we abbreviate $x^i\1+\R$  to $x^i\1$ and refer to the basis $\{x^i\1| 0\le i\le e\}$ of $L(\lambda)$ as {\em standard}. Let $\widehat{H}$ denote the set of isomorphism classes of simple $H$- modules and $\widehat{H}_e$ be the subset of modules $L(\lambda)$ with $e(\lambda)=e$. We comment on the distribution of simple modules among classes $\widehat{H}_e$.
\begin{Prop}\label{classesofsimples} Suppose $\k$ contains a root of unity of order equal to the exponent of $G$.
\begin{enumerate}
\item[(1)] Every $\widehat{H}_e$ is nonempty for all $0\le e\le n-1$ iff $|ab|$ is divisible by $n$,
\item[(2)] The cardinalities of the nonempty $\widehat{H}_e$ are equal for all $e<n-1$,
\item[(3)] Let $n'= n/2$ whenever n is even. If $|ab|$ is not divisible by $n$, then $\widehat{H}_e$ is nonempty for even $e$ only. $\widehat{H}_{n-1}$ is nonempty iff $|ab|\ne n'$.
\end{enumerate}\end{Prop}
\pf (1) $\widehat{H}_e\ne\emptyset$ for $0\le e<n-1$ iff the equation
$$\lambda(ab)=q^e$$
is solvable in $\G$. Further, $\widehat{H}_{n-1}\ne\emptyset$ iff either $\lambda(ab)=q^{n-1}$ is solvable or $\exists\lambda\in\G:\lambda(ab)\notin\gen q$. Therefore, if $\widehat{H}_e\ne\emptyset$ for all $e$ then $\lambda(ab)=q$ is solvable, hence $|ab|$ is divisible by $n$.

Conversely, if $|ab|$ is divisible by $n$, there is a homomorphism $\check{\lambda}:\gen{ab}\to\gen q,\ \check{\lambda}(ab)=q^e$ for every $0\le e\le n-1$. Let $U$ be the subgroup of $\k^{\bullet}$ generated by an element of order equal to the exponent $E$ of $G$. $U$ is injective in the category of abelian groups of exponent $E$. Therefore $\check{\lambda}$ can be lifted to $\lambda:G\to\k^{\bullet}$, which completes the proof of (1).

(2) The cardinality of every nonempty $\widehat{H}_e,\ 0\le e< n-1$ equals the order of $(ab)^{\perp}:=\{\phi\in\G|\phi(ab)=1\}$.

(3) From $\chi(ab)=q^{-2}$ we see that $|ab|$ is divisible by $|q^2|$. The latter is $n$ iff $n$ is odd, and $n'$, otherwise. Therefore in this case $|ab|=n'm$ for some integer $m\ge 1$. If $m=1$, then every homomorphism $\gen {ab}\to\gen q$ sends $ab$ to $q^{2i}$ for some $0\le i<n'$ and those are the only possibilities. If $|ab|\ne n'$, there are $\lambda\in\G$ such that $\lambda(ab)\notin\gen q$.\qed
\begin{Prop}\label{projective} $L(\lambda)$ is projective iff $e(\lambda)=n-1$.
\end{Prop}
\pf In one direction, if $e(\lambda)< n-1$, then by the preceeding proposition $\Z$ is a nonsplit extension of $\R$ by $L(\lambda)$, thus $L(\lambda)$ is nonprojective.

Conversely, let $e(\lambda)=n-1$. For every $\lambda\in\G$ define a primitive idempotent of $\k G$ by
$$\e=|G|^{-1}\sum_{g\in G}\lambda(g^{-1})g.$$
Now pick the element $w:=\e y^{n-1}x^{n-1}\in H$. Clearly $gw=\lambda(g)w$ which means that $w$ has weight $\lambda$. Further we note the identity
\begin{align*}y\e&=|G|^{-1}(\sum_{g\in G}\lambda(g^{-1})\chi(g)g)y= |G|^{-1}(\sum_{g\in G}\lambda\chi^{-1}(g)g)y\\
&=e_{\lambda\chi^{-1}}y\end{align*} 
where we used the equality $yg=\chi(g)gy$. We see readily that $w$ is a $\lambda$- primitive element, and therefore $Hw\simeq \Z$, as $\Z$ is simple. Thus to complete the proof it suffices to show that $w^2=cw$ for some $0\ne c\in\k$. To this end we employ a formula \cite[1.3]{KR} for commuting powers of $y$ and $x$. To adjust \cite[Lemma 1.3]{KR} to the present situation we replace $q$ by $q^{-1}$ and set the coefficient $\lambda=-q^{-1}$. Then we get           $$y^jx^k-q^{-jk}x^ky^j=\sum_{i=1}^{\text{min}(j,k)}x^{k-i}f_i^{j,k}y^{j-i}$$                
where                
$$f_i^{j,k}=(-q^{-1})^i{\binom ji}_{q^{-1}}{\binom ki}_{q^{-1}}(i)_{q^{-1}}!q^{-(k-i)(j-i)}\prod_{m=1}^i(q^{m+i-j-k}ab-1).$$
We calculate $w^2$ keeping in mind that $yw=0$, namely
\begin{align*}w^2&=e_{\lambda}y^{n-1}x^{n-1}w=e_{\lambda}f_{n-1}^{n-1,n-1}w\\
&=(-q^{-1})^i(n-1)_{q^{-1}}!\prod_{m=1}^{n-1}(q^{m-n+1}\lambda(ab)-1)w.\end{align*}
As $\lambda(ab)=q^{n-1}$ or $\lambda(ab)\not\in\gen q$ at all, $\prod_{m=1}^{n-1}(q^{m-n+1}\lambda(ab)-1)\ne 0$ which completes the proof.\qed

\section{Trace of the Antipode}\label{trace}

We begin with a discussion of simple finite-dimensional modules for an algebra $H$. Let $V$ be one such module. Denote by  $P_V=\text{ann}_H V$ the primitive ideal associated to $V$. The next observation is well-known.
\begin{Lem}\label{standardfact} Let $H$ be a Hopf algebra with invertible antipode $S$, $V,V'$ two simple $H$- modules. $V^*$ is isomorphic to $V'$ iff $S(P_{V'})=P_V$.
\end{Lem}
\pf If simple $H$- modules $V$ and $W$ are isomorphic, then $\text{ann}_HV=\text{ann}_HW$. Conversely, if $P_V=P_W$, then $V$ and $W$ are simple modules for 
a simple algebra $H/P_V$, hence they are isomorphic.

Next we compute $\text{ann}_HV^*$. For every $h\in H$, $HV^*=0$ iff\newline $f(S(h).v)=0$ holds for all $f\in V^*$ and $v\in V$. Thus $S(h).v=0$ for all $v\in V$, whence $S(P_{V^*})\subset P_V$. Since $P_{V^*}$ and $P_V$ are maximal ideals of $H$, so are $S(P_{V^*})$ and $P_V$, and thus they are equal. We see that $V^*\simeq V'$ iff $S^{-1}(P_V )=P_{V'}$.\qed

The above Lemma says that for a cofinite primitive ideal $P$ of $H$, $P$ is fixed by $S$ iff $V$ is self-dual. When this is the case, $S$ induces a linear transformation, denoted by $S_V$, on $H/P\cong \text{End}_{\Delta}(V)$, where $\Delta$ is the centralizer of $H$ in $\text{End}_{\k}(V)$. We are interested in computing $\text{Tr}\ (S_V)$ for all simple self-dual $H=H(\D)$- modules. We move on to a description of these modules for this $H$. 
\begin{Lem}\label{dual} For every $\lambda\in\G$, $L(\lambda)^*\simeq L(\lambda^{-1}\chi^{-e})$ where $e=e(\lambda)$.
\end{Lem}
\pf We pick the functional $f_0\in L(\lambda)^*$ dual to $x^e\1$, i.e. given $0\leq i\leq e,$ $f_0(x^i\1)=\delta_{ie}$. From 
$$g.f_0(x^i\1)=f_0(g^{-1}.x^i\1)=\chi^i(g^{-1})\lambda(g^{-1})\delta_{ie}$$
we see that $f_0$ has weight $\lambda^{-1}\chi^{-e}$. $f_0$ is also primitive by the following calculation 
$$y.f_0(x^i\1)=f_0(S(y).x^i\1)=-f_0(b^{-1}y.x^i\1)=0,$$
because by \eqref{simplecommutator} $y.x^i\1$ is a scalar multiple of $x^{i-1}\1$.\qed

\begin{Cor}\label{selfdual}\begin{itemize}
\item[(1)] $L(\lambda)$ is self-dual iff $\lambda^2=\chi^{-e}$,
\item[(2)] If $|G|$ is odd, every class $\widehat{D}_e$ has a unique self-dual module.
\end{itemize}
\end{Cor}
\pf Part (1) follows immediately from the previous Lemma.

(2) The order of $\G$ is also odd, hence the mapping $\lambda\to \lambda^2$ is an isomorphism. Therefore every equation $\lambda^2=\mu$ has a unique solution.\qed

We digress to discuss a relationship between various spaces of homomorphisms. Let $H$ be a Hopf algebra and $V$ an $H$- module. A bilinear map $\gen{\cdot,\cdot}:V\o V\to\k$ is called $S$- adjoint if $h$ and $S(h)$ are adjoint operators relative to the form for every $h$, i.e.
\begin{equation}\label{Sadjoint}\gen{hv,w}=\gen{v,S(h)w}.\end{equation}
We denote by $\text{Hom}_{\k}^S(V\o V,\k)$ the space of all $S$- adjoint bilinear maps on $V$. Every $f\in\text{Hom}_{\k}(V,V^*)$ gives rise to a bilinear form by
$$\gen{\cdot,\cdot}:V\o V\to\k,\;\gen{v,w}=f(v)(w)$$
and, conversely, a form $\gen{v,w}$ defines an $f$ by that formula. Moreover, $f$ is $H$- linear, i.e. $f(h.v)=hf(v)$, iff the associated form is $S$- adjoint. We note that $\text{Hom}_{\k}^S(V\o V,\k)$ contains the subspace $\text{Hom}_{\k}(V\o V,\k)^H$ of all $H$- invariant forms, and the latter is isomorphic to $\text{Hom}_{\k}(V,V^*)^H$, provided $S$ is invertible. When $H$ is cocommutative all four spaces can be identified.

Fix a basis of $V$ and denote by $[v]$ the column vector representing $v$ and by $[h]$ the matrix of $h$ relative to the basis. Let $A=(a_{ij})$ be the matrix of $\gen{\cdot,\cdot}$ in that basis, namely
$$\gen{v,w}=[v]^TA[w].$$ 
\begin{Lem}\label{formulafortrace} Suppose $V$ is a simple, absolutely irreducible, self-dual $H$- module. Let $A$ be the matrix of the bilinear form associated to an isomorphism $V\to V^*$. Then
\begin{equation}\label{traceformula} \text{Tr}\ (S_V)=\text{Tr}\ (A^{-1}A^T)\end{equation}
\end{Lem}
\pf Equality \eqref{Sadjoint} says in the matrix form
$$[v]^T[h]^TA[w]=[v]^TAS([h])[w]$$
which is equivalent to $[h]^TA=AS([h])$. Set $X=[h]$. By our assumptions on $V$, $X$ runs over all $d\times d$- matrices, where $d=\text{dim}\ V$. It follows that the action of $S$ is described by 
$$S(X)=A^{-1}X^TA.$$
Denote by $\hat{a}_{ij}$ the $ij$- entry of $A^{-1}$. It is elementary to see that
$$\text{Tr}\ (S_V)=\sum_{ij}\hat{a}_{ij}a_{ij}=\text{Tr}(A^{-1}A^T).$$
An alternative proof of the Lemma can be found in \cite{MS}.\qed

We return to a self-dual $H$- module $L(\lambda)$. By Lemma \ref{dual} $\psi:L(\lambda)\to L(\lambda)^*,\;\psi(x^i\1)=x^i.f_0$ is an $H$- isomorphism.
\begin{Prop}\label{matrix} Let $\gen{\cdot,\cdot}:L(\lambda)\o L(\lambda)\to\k$ be the $S$-adjoint form induced by $\psi$ and $A=(a_{ij})$ be its matrix relative to the canonical basis of $L(\lambda)$. Then
$$a_{ie-i}=(-1)^iq^{-\binom i2+ei}\lambda(a)^{-i}.$$
\end{Prop}
\pf By definition of $A$
$$a_{ij}=\gen{x^i\1,x^j\1}=x^i.f_0(x^j\1)=f_0(S(x^i)x^j\1).$$
To continue the calculation we use the well-known equality $S(x^i)=(-a^{-1}x)^i=(-1)^iq^{-\binom i2}a^{-i}x^i$:
\begin{align*}\label{entry} a_{ij}&=(-1)^iq^{-\binom i2}f_0(a^{-i}x^{i+j}\1)\\
&=(-1)^iq^{-\binom i2}\chi^{i+j}(a^{-i})\lambda(a^{-i})f_0(x^{i+j}\1)\\
&=(-1)^iq^{-\binom i2 +ei}\lambda(a)^{-i}\delta_{i+j,e}\end{align*}
as $\chi(a)=q^{-1}$.\qed

We proceed to the main theorem of this note.
\begin{Thm}\label{main} For every $\lambda\in\G$ with the associated integer $e=e(\lambda),$
$$\text{Tr}\ (S_{L(\lambda)})=(-1)^e\lambda(a)^eq^{-\binom {e+1}2}(e+1)_q.$$
\end{Thm}
\pf We find it convenient to denote an $(e+1)\times(e+1)$ matrix $A$ whose only nonzero entries are $a_{ie-i},i=0,\ldots,e$ by 

\noindent$\text{diag}'(a_{0e},\cdots,a_{ie-i},\cdots,a_{e0})$. It can be seen easily that both $A^{-1}$ and $A^T$ are of that same form, explicitly one has
$$A^{-1}=\text{diag}'(\hat{a}_{0e},\ldots,\hat{a}_{ie-i},\ldots,\hat{a}_{e0})\;\text{and}\;A^T=\text{diag}'(a'_{0e},\ldots,a'_{ie-i},\ldots,a'_{e0})$$
with $\hat{a}_{ie-i}=a^{-1}_{e-ii}$ and $a'_{ie-i}=a_{e-ii}$. Their product is given by
$$A^{-1}A^T=\text{diag}\ (\hat{a}_{0e}a'_{e0},\ldots,\hat{a}_{ie-i}a'_{e-ii},\ldots,\hat{a}_{e0}a'_{0e})$$
and therefore
\begin{equation}\label{matrixtrace}\text{Tr}\ (A^{-1}A^T)=\sum_{i=0}^e\hat{a}_{ie-i}a'_{e-ii}=\sum_{i=0}^ea^{-1}_{e-ii}a_{ie-i}\end{equation}
From Proposition \ref{matrix} we have 
\begin{align*}a^{-1}_{e-ii}a_{ie-i}&=(-1)^{e-i}q^{\binom{e-i}2-e(e-i)}\lambda(a)^{e-i}(-1)^iq^{-\binom i2+ei}\lambda(a)^{-i}\\&=(-1)^e\lambda(a)^eq^{s_i}\lambda(a)^{-2i},\end{align*}
where $s_i=\binom{e-i}2-e(e-i)-\binom i2+ei$. Since $L(\lambda)$ is self-dual, $\lambda^2=\chi^{-e}$, which gives $\lambda^2(a^{-i})=\chi(a)^{ei}=q^{-ei}$. From this we derive $a^{-1}_{e-ii}a_{ie-i}=(-1)^e\lambda(a)^eq^{s'_i}$ with $s'_i=\binom {e-i}2-e(e-i)-\binom i2$. An elementary calculation gives $s'_i=i-\binom{e+1}2$. Summing up we get
$$\text{Tr}\ (S_{L(\lambda)})=(-1)^e\lambda(a)^eq^{-\binom{e+1}2}\sum_{i=0}^eq^i=(-1)^e\lambda(a)^eq^{-\binom{e+1}2}(e+1)_q.$$\qed

Combining Proposition \ref{projective} with the last theorem we obtain
\begin{Cor}\label{traceprojective} For a self-dual module $L$,$\text{Tr}\ (S_L)=0$ iff $L$ is projective.
\end{Cor}

\section{Special Cases}\label{special} The trace invariant has been computed for modules $V$ of certain Hopf algebras. We want to show that Theorem \ref{main} contains these previously obtained formulas for $Tr(S_V)$ .

1. First we take up the small quantum group for $\mathfrak{sl}_2$, $H=\bf{u}_q(\mathfrak{sl}_2)$. The associated datum is $\D=(G,K,K,\chi,1)$ where $G=\gen{K|K^n=1},\newline \chi(K)=r^2$ with $r$ a root of $1$ of order $n$, $n$ an odd integer. The parameter $q$ of the general definition of lifting in Section \ref {liftings} is $q=\chi(K^{-1})=r^{-2}$. For a character $\lambda:G\to\k^{\bullet}$ the integer $e=e(\lambda)$ is defined by
$$\lambda(K^2)=q^e=r^{-2e}.$$
We see that $\lambda(K)=r^{-e}$ as the other solution $\lambda(K)=-r^{-e}$ violates the equality $\lambda(K)^n=1$. Set $L(e)=L(\lambda)$, for the isomorphism class of every simple module is determined by $e=e(\lambda)$. By Theorem \ref{main}, substituting $q$ by $r^{-2}$, we have
$$\text{Tr}\ (S_{L(e)})=(-1)^er^{-e^2}r^{2\binom{e+1}2}r^{-e}[e+1]_r=(-1)^e[e+1]_r,$$
where $[m]_r:=\displaystyle\frac{r^m-r^{-m}}{r-r^{-1}}$.
This agrees with $\mu(V_e)=Tr(S_{V_e})$ as determined in \cite{J} where, in the notation of that paper,  $V_e$ is the self-dual $(e+1)$-dimensional $\bf{u}_q(\mathfrak{sl}_2)$-module that corresponds to our $L(e)$.

2. The next case is the quantum double of a Borel-like subalgebra of $\bf{u}_q(\mathfrak{sl}_2)$. Let $H=\bf{u}_q^{\ge 0}$ be the sub-Hopf algebra generated by $K$ and $E$. Note that $H$ is a Taft algebra with parameters $r^2$ and $n=|r|$ as $n$ is assumed to be an odd integer. $H$ is also a rank $1$ Hopf algebra \cite{KR1} with the group of grouplikes $G=\gen K$ and the datum $(G,K,\chi)$ where $\chi:G\to\k^{\bullet},\ \chi(K)=r^2$. Set $D=D(H)$ for the Drinfel'd double of $H$. We give a description of $D$ following \cite[Prop. 5]{KR1}. Let $\G=\text{Hom}\ (G,\k)$ and $\Gamma=G\times\G$. Pick $\phi:\Gamma\to\k^{\bullet}$ defined by $\phi(g\gamma)=\chi(g)\gamma(K^{-1}),g\in G,\gamma\in\G$. In keeping with conventions of \cite{KR1} we write $\phi=\chi\widehat{K}^{-1}$, the point is that every $g\in G$ induces a homomorphism $\hat{g}:\G\to\k, \hat{g}(\gamma)=\gamma(g)$. Pick  a functional $F: H\to\k$ defined by the equality
$$F(E^ig)=\delta_{i1}\;\text{for all}\;g\in G.$$
By \cite[Prop. 5]{KR1} $D$ is generated by $\Gamma,E,F$ subject to the relations
\begin{align}K^n&=1=\chi^n\\ (g\gamma)E&=\phi(g\gamma)E(g\gamma)\;\text{for all}\; g\gamma\in\Gamma\\
(g\gamma)F&=\phi^{-1}(g\gamma)F(g\gamma)\;\text{for all}\; g\gamma\in\Gamma\\
E^n&=0=F^n\\ &EF-FE=\chi-K\label{commutator}\\
\Delta(E)&=K\o E+E\o 1\\\Delta(F)&=1\o F+F\o\chi\label{comultiplication}\end{align}
$D(H)$ is actually a lifting of a quantum plane. To see this we introduce a new generator $F':=-F\chi^{-1}$. One can check easily the equalities
$$EF'-qF'E=K\chi^{-1}-1\;\text{with}\;q=\phi^{-1}(K)=r^{-2},\;\text{and}$$
$$\Delta(F')=\chi^{-1}\o F'+F'\o 1$$
which replace \eqref{commutator} and \eqref{comultiplication}, respectively. Comparing these defining relations of $D$ with the defining relations of a lifting of a quantum plane we see that $D$ is a lifting associated to the datum $(\Gamma,K,\chi^{-1}\phi,1)$.

For a character $\lambda:\Gamma\to\k^{\bullet}$ the integer $e=e(\lambda)$ is defined by 
\begin{equation}\label{specialexponent}\lambda(K\chi^{-1})=q^e=r^{-2e}.\end{equation}
because the order of $K\chi^{-1}$ is $n$, and so is the order of $q$. We want to determine the weights of self-dual modules. That is to say, by Corollary \ref{selfdual}(1) we seek solutions to the equation $\lambda^2=\phi^{-e}$. Since $n$ is odd $\gen{K\chi^{-1}}\cap\gen{K\chi}=1$, hence $\Gamma=\gen{K\chi^{-1}}\times\gen{K\chi}$ is the direct product of those subgroups. Let $\{\alpha,\beta\}$ be the dual basis of $\widehat{\Gamma}$. Namely, $\alpha,\beta\in\widehat{\Gamma}$ are such that
\begin{align*}\alpha(K\chi^{-1})&=r^2,\;\alpha(K\chi)=1\\
\beta(K\chi^{-1})&=1,\;\beta(K\chi)=r^2\end{align*}
Right from the definition of $\phi$ we have $\phi(K\chi)=1$ and $\phi(K\chi^{-1})=r^4$ which gives $\phi=\alpha^2$. As for $\lambda$, assuming $\lambda(K\chi)=r^{2j}$ we have by \eqref{specialexponent}, $\lambda=\alpha^{-e}\beta^j$. The condition of Corollary \ref{selfdual}(1) gives $\alpha^{-2e}\beta^{2j}=\alpha^{-2e}$, whence $j=0$. It follows that the self-dual modules are of the form $L(e):=L(\alpha^e),0\le e\le n-1$. To determine the scalar $\alpha^e(K^2)$ that appears in the trace formula we compute $\alpha^e(K^2)=\alpha(K\chi K\chi^{-1})=r^{2e}$, hence $\alpha(K)=r^e$. Therefore by Theorem \ref{main}, setting $q=r^{-2}$ there, we have
\begin{equation*}\text{Tr}\ (S_{L(e)})=(-1)^er^{e^2}r^{-2\binom{e+1}2}(e+1)_{r^{-2}}=(-1)^e[e+1]_r.\end{equation*}

\end{document}